\documentclass[a4paper,11pt]{article}

\usepackage{amsfonts} 
\usepackage{amssymb} 
\usepackage{amsmath,amsthm} 

\newtheorem{theoreme}{Theorem}[section] 
\newtheorem{lemma}[theoreme]{Lemma} 
\newtheorem{prop}[theoreme]{Proposition} 
\newtheorem{cor}[theoreme]{Corollary} 
\newtheorem{definition}[theoreme]{Definition} 
 
 \theoremstyle{remark}

\def \Rk {\ {\bf Remark.} } 
 
\def \sm {\setminus } 
\newcommand{\be}{\begin{enumerate}}  \newcommand{\ee}{\end{enumerate}} 
\newcommand{\bi}{\begin{itemize}}  \newcommand{\ei}{\end{itemize}} 
\newcommand{\bd}{\begin{description}}  \newcommand{\ed}{\end{description}}

\newcommand{\comment}[1]{}

\def \R {\mathbb{R}}   
 \def \T {\mathbb{T}} \def \H {\mathbb{H}}

\DeclareMathOperator{\vol}{vol}

\DeclareMathOperator{\diam}{diam}
\DeclareMathOperator{\cotanh}{cotanh}

\numberwithin{equation}{section}  
\renewcommand{\phi}{\varphi} 
\renewcommand{\epsilon}{\varepsilon} 
 
\title{On the fundamental group of some open manifolds}
\author{Nader Yeganefar}
\date{\today} 

\begin{document} 
\maketitle
\begin{abstract}
We study fundamental groups of non compact Riemannian manifolds. We find conditions which ensure that the fundamental group is trivial, finite or finitely generated.
\end{abstract}


\section{Introduction}
 
Let $M$ be a complete non compact Riemannian manifold. It is a classical theme in Riemannian geometry to find geometric conditions which ensure some finiteness results for the topology of $M$. For example, what conditions would imply that $M$ has finite topological type (i.e. is homeomorphic to the interior of a compact manifold with boundary)? Answers to this question are numerous and can be easily found in the literature, but we won't give any reference here. 

In this note, we are merely interested in the fundamental group of $M$. As $M$ is not compact, the basic problem is to know if its fundamental group is finitely generated or not, and then to know if it is finite or even trivial. Here, we will deal mainly with two situations.

First, we assume that $M$ has sectional curvature $K$ bounded below by a negative constant, i.e. $K\geq -1$. Denote by $V_{\H^n}(r)$ the volume of a ball of radius $r$ in hyperbolic space $\H^n$. Then 
$$V_{\H^n}(r)=\omega _{n-1}\int_0^r\sinh ^{n-1}{(t)}\, dt,$$ 
where $\omega _{n-1}$ is the volume of the unit sphere $S^{n-1}\subset \R ^n.$ The Bishop-Gromov theorem \cite{Gr2} asserts that for all $p\in M$, the function $r\mapsto \vol{(B(p,r))}/V_{\H^n}(r)$ is decreasing, where  $B(p,r)$ is a ball of radius $r$ around $p$ in $M$. (Of course, a lower bound on Ricci curvature is enough for this.) Therefore, the quantity
$$v_p:=\lim _{r\to \infty}{\frac{\vol{(B(p,r))}}{V_{\H ^n}(r)}}\in[0,1]$$ 
is well defined. We also set
$$v(M)=\inf _{p\in M}{v_p}.$$
Manifolds with $v(M)>0$ are said to have large volume growth and were studied by Xia \cite{X2}, who  found additional conditions which imply that $M$ has finite topological type or is diffeomorphic to $\R ^n$. Here, our first result is 
\begin{theoreme}\label{length}
Let $M^n$ be a complete $n-$dimensional non compact Riemannian manifold, and fix a point $p\in M$. Assume that the sectional curvature $K$ is bounded below $K\geq -1$. For any real number $L>0$, there exists a constant $v=v(L,n)\in (0,1)$ such that if
$$\lim _{r\to \infty}{\frac{\vol{B(p,r)}}{V_{\H ^n}(r)}}\geq 1-v$$
and $M$ is not simply connected, then the length of the shortest homotopically non trivial geodesic loop based at $p$ is bigger than $L$.
\end{theoreme}
For $r>0$, set
$$\delta (r)=\frac{\sqrt{\cosh{(r)}}}{\cosh{(r/2)}}.$$
Combining our theorem with techniques developed by Xia in \cite{X2}, we get
\begin{cor}\label{Xia}
Let $(M^n,g)$ be a complete $n-$dimensional non compact Riemannian manifold, and fix a point $p\in M$. Assume that the sectional curvature $K$ is bounded below $K\geq -1$. Given a real number $L>0$, there exists $v=v(L,n)\in (0,1)$ such that if
\begin{itemize}
\item[i)] $$\lim _{r\to \infty}{\frac{\vol{B(p,r)}}{V_{\H ^n}(r)}}\geq 1-v,$$
\item[ii)] $M$ has large volume growth, i.e. $v(M):=\inf _{q\in M}{v_q}>0$,
\item[iii)] and for all $r\geq L$,
$$\frac{\vol{B(p,2r)}-v_p V_{\H ^n}(2r)}{v(M)\omega _{n-1}}<\int _0^{\cosh ^{-1}{(\delta(2r))}} \sinh ^{n-1}{(t)} \, dt,$$
\end{itemize}
then $M$ is simply connected.
\end{cor}

Next we consider manifolds without curvature bounds but we impose some growth conditions at infinity. More specifically, if $M$ is any complete manifold and $p\in M$ is any point, we can study the asymptotic behavior of the diameter of geodesic spheres around $p$, as the radius goes to infinity. Let $\diam{\partial B(p,r)}$ denote the diameter of a sphere of radius $r$ around $p$, measured with respect to the distance of $M$. By the triangle inequality, we have always 
$$\diam{\partial B(p,r)}\leq 2r.$$
For example, if there is a line in $M$ passing through $p$, then the diameter of a sphere of radius $r$ around $p$ is exactly $2r$. On the other hand, in \cite{Sh} Z. Shen considered manifolds satisfying  
\begin{equation}\label{growth}
\limsup _{r\to \infty}{\frac{\diam{\partial B(p,r)}}{r}}<1.
\end{equation}
Shen showed that a complete manifold satisfying this smallness condition at infinity is proper (i.e. its Busemann functions are proper). In another direction, C. Sormani \cite{S} proved that for each integer $n$, there exists an explicit small constant $\epsilon _n>0$ such that if $M$ is a complete $n-$dimensional manifold with nonnegative Ricci curvature, which for some point $p$ satisfies
$$\limsup _{r\to \infty}{\frac{\diam{\partial B(p,r)}}{r}}<\epsilon _n,$$
then it has finitely generated fundamental group.

Our next result implies that if the universal cover of a complete non compact manifold $M$ is not too big at infinity, then $M$ has finite fundamental group.
\begin{theoreme}\label{finite}
Let $\tilde{M}$ be a complete non compact Riemannian manifold and let $\tilde{p}\in \tilde{M}$ be a point. Let $G$ be a discrete group of isometries acting freely on $\tilde{M}$. If $\tilde{M}$ has small diameter growth
$$\limsup _{r\to \infty}{\frac{\diam{\partial B(\tilde{p},r)}}{r}}<1,$$
then either $\tilde{M}/G$ is compact or $G$ is finite.
\end{theoreme}
The following question is implicit in \cite[p. 291]{Sh} (after Corollary 1): is it true that any complete non compact manifold of positive Ricci curvature has small diameter growth in the sense of (\ref{growth})? Some partial (positive) answer to this question is being worked out by Cao and Wang \cite{CW}, but we will see in a moment that the answer is in general negative. Namely, it is a consequence of the work of Belegradek and Wei \cite{BW} that if $k$ is a given integer, then for a sufficiently large integer $d$, there is a metric of positive Ricci curvature on $\T ^k\times \R ^d$, where $\T ^k$ is a flat torus. Hence there is also a metric of positive Ricci curvature on the universal cover $\R ^{k+d}$. As the fundamental group of $\T ^k\times \R ^d$ is infinite, a corollary of Theorem \ref{finite} is
\begin{cor}
For $n$ sufficiently large, there is a metric of positive Ricci curvature on $\R ^n$ which does not have small diameter growth (\ref{growth}).
\end{cor}

Now consider any Riemannian manifold $M$, and let $(p,\gamma _1$, $\gamma _2)$ be a hinge in $M$, i.e. $\gamma _1$ and $\gamma _2$ are two minimizing geodesics such that $p=\gamma _1(0)=\gamma _2(0)$. Assume that $\gamma _1$ and $\gamma _2$ have the same length, to be denoted by $L>0$.
\begin{definition}\label{hinge}
We say that $M$ has the thick hinges property at $p$ if there exists a constant $\theta _0>0$ such that if $(p,\gamma _1$, $\gamma _2)$ is any hinge as described above which satisfies $d(\gamma _1(L), \gamma _2(L))\geq L$, then the angle of the hinge at $p$ is greater than or equal to $\theta _0$.
\end{definition}  
For example, if $M$ has nonnegative sectional curvature then by the Toponogov comparison theorem $M$ has the thick hinges property at any point $p$ (and we can choose $\theta _0=\pi /3$.) It would be interesting to find other less restrictive geometric conditions which imply this property. Our final result is
\begin{theoreme}\label{fg}
Let $M$ be a complete non compact Riemannian manifold and let $p\in M$ be a point. Consider the universal Riemannian cover $\tilde{M}$ and choose a lift $\tilde{p}\in \tilde{M}$ of $p$. If $\tilde{M}$ has the thick hinges property at $\tilde{p}$, the fundamental group $\pi _1(M,p)$ is finitely generated. 
\end{theoreme}

The paper is organized as follows. In section \ref{fund}, we gather some well-known facts about fundamental groups of Riemannian manifolds. In section 3, we prove Theorem \ref{length} and Corollary \ref{Xia}. In the last section, we prove Theorem \ref{finite}.\\

\noindent {\bf Acknowledgements.} I would like to thank Erwann Aubry for useful comments and hints. In particular, he pointed out to me this version of Theorem \ref{finite} which improves a previous version, and brought to my attention the work of Belegradek and Wei. Thanks are also due to Zhongmin Shen and especially to Jianguo Cao for communicating me his current work with Jiaping Wang \cite{CW}. Finally, I thank Gilles Carron for his encouragements and comments.


\section{Fundamental group and generators}\label{fund}
In this section, we review some standard facts about fundamental groups of Riemannian manifolds (see e.g. \cite{Gr}, \cite{GP} and \cite{S}). Although the proofs are easy, we reproduce them here for completeness. Let $M$ be any complete Riemannian manifold, and fix a point $p\in M$. We assume that $M$ is not simply connected. If $C\in \pi _1(M,p)$ is a homotopy class of loops based at $p$, its length $L(C)$ is by definition the infimum of the lengths of loops in $C$. By an elementary argument, there exists a geodesic loop $c\in C$ whose length is precisely $L(C)$; we will refer to such a loop $c$ as a \emph{minimal representative geodesic loop}. Note that $c$ is not necessarily smooth at $p$. 
\begin{lemma}\label{fini}
For every real number $r>0$, the set $\{ L(C): C\in \pi _1(M,p),\\ L(C)\leq r \}$ is finite.
\end{lemma}
\begin{proof}
Assume on the contrary that for some $r>0$, we have an infinite sequence $C_i$ of mutually distinct elements in $\pi _1(M,p)$ with $L(C_i)\leq r$. For each $i$, choose a minimal representative geodesic loop $c_i\in C_i$. We may suppose that each $c_i$ has unit speed. Therefore, by Ascoli theorem, there is a subsequence of $\{ c_i\}$ which converges uniformly. In particular, $c_i$ and $c_j$ are as close as we wish, provided $i\neq j$ are suitable large indices. But two sufficiently close loops are homotopic, and hence $c_i$ and $c_j$ are in the same homotopy class for some $i\neq j$. This is a contradiction.
\end{proof}

\begin{definition}
We say that a non trivial element $C$ in $\pi _1(M,p)$ is irreducible when each decomposition $C=C_1C_2$ is such that $L(C_1)\geq L(C)$ or $L(C_2)\geq L(C)$.
\end{definition}
For example, if $c$ is a non contractible geodesic loop based at $p$ whose length is minimal, then it represents an irreducible classe. 
An interesting geometric property of irreducible classes is the following:
\begin{prop}\label{minimal}
Let $C\in \pi _1(M,p)$ be an irreducible class, with minimal representative geodesic loop $c$. Then $c$ is a minimizing geodesic on $[0,L(C)/2]$ and on $[L(C)/2,L(C)]$.
\end{prop}
\begin{proof}
For any curve $\gamma \colon [0,L]\to M$, denote by $-\gamma$ the curve defined by $-\gamma (t)=\gamma (L-t)$. Now we argue by contradiction and assume that there exists a minimizing geodesic $c'$ from $p$ to $c(L(C)/2)$ which has length less than $L(C)/2$.  Then the loops $c_1$ and $c_2$, obtained by composing $c\vert _{[0,L(C)/2]}$ with $c'$ and respectively $-c'$ with $c\vert _{[L(C)/2,L(C)]}$ have length less that $L(C)$, and $c$ can be homotoped to the product $c_1c_2$. Thus $C$ is not irreducible, which is a contradiction.
\end{proof}
It follows from this proposition that the midpoint $c(L(C)/2)$ is in the cut locus of $p$, because there are two minimizing geodesics running from $p$ to this point. We have actually a bit more, as we will explain now. Consider the function $d(p,.)$ defined on $M$. Although it is not a smooth function, there is a notion of critical point for $d(p,.)$ which generalizes the usual notion. This has been introduced by Grove and Shiohama (see the nice surveys of Cheeger \cite{C} and Grove \cite{Gro}). More specifically, a point $m\in M$ is called \emph{critical} for $d(p,.)$ if for each tangent vector $u$ at $m$, there exists a minimizing geodesic from $m$ to $p$ whose initial velocity makes an angle $\leq \pi /2$ with $u$. This definition allows us to develop a kind of Morse theory for $d(p,.)$. 
Now in the situation described in the above proposition, there are two minimizing geodesics running from $c(L(C)/2)$ to $p$, and their angle at $c(L(C)/2)$ is $\pi$. Therefore, we easily get
\begin{cor}\label{critic}
Let $C\in \pi _1(M,p)$ be an irreducible class, with minimal representative geodesic loop $c$. Then $c(L(C)/2)$ is a critical point for $d(p,.)$.
\end{cor}
\Rk There is no general definition of index of a critical point of $d(p,.)$, as in usual Morse theory. If such a definition existed, we would like $c(L(C)/2)$ to be a critical point of index $1$.\\

Our interest in irreducible classes comes from the following fact:
\begin{prop}\label{generators}
$\pi _1(M,p)$ is generated by its irreducible classes. 
\end{prop}
\begin{proof}
Set $G=\pi _1(M,p)$. We observe first the following fact: let $H$ be a subgroup of $G$ such that $G\sm H\neq \emptyset$, and let $C$ be an element $G\sm H$ which minimizes length, i.e.
\begin{equation}\label{aaaa}
\forall C'\in G\sm H,\,\, L(C)\leq L(C').
\end{equation}
Then $C$ is irreducible. Namely, assume on the contrary that we can write $C=C'C''$, with $L(C')<L(C)$ and $L(C'')<L(C)$. Then by \eqref{aaaa}, $C'$ and $C''$ are elements of $H$, so that $C$ itself is in $H$. This contradicts the definition of $C$.

Now, we define a sequence of generators $\{ C_i\}$ of $G$ as follows. Choose an element $C_1\in G$ such that it minimizes length, i.e.
$$\forall C\in G\sm \{1\},\,\, L(C_1)\leq L(C).$$
Let $H_1$ be the subgroup generated by $C_1$. Define $C_i$ inductively as follows. Assume that we have elements $C_1,\, C_2,\, \ldots ,\, C_{i-1}$ which generate a subgroup $H_{i-1}$. If $G\sm H_{i-1}$ is not empty, choose an element $C_i$ in $G\sm H_{i-1}$ such that
$$\forall C\in G\sm H_{i-1},\,\, L(C_i)\leq L(C).$$
By the discussion above, each $C_i$ is irreducible. Moreover, we have then two possibilities. Either the sequence $L(C_i)$ is bounded, or it goes to infinity. In the first case, there is only a finite number of $C_i$'s by lemma \ref{fini} and we must have $G=H_i$ for some $i$. In the second case, for each element $C$ of $\pi _1$, there is an $i$ such that $L(C)<L(C_i)$, so that $C\in H_i$; hence the $C_i$'s generate $G$.
\end{proof}
\noindent\Rk These irreducible generators are called "halfway generators" by Sormani, see \cite[Lemma 5, Definition 6]{S}.

\section{Proofs of Theorem \ref{length} and Corollary \ref{Xia}}

\begin{proof}[Proof of Theorem \ref{length}]
We will prove that there exists a constant $v=v(L,n)\in (0,1)$ such that if there is a non trivial element $C$ in $\pi _1(M,p)$ of length $R:=L(C)\leq L$, then
$$\lim _{r\to \infty}{\frac{\vol{B(p,r)}}{V_{\H ^n}(r)}}\leq 1-v.$$
Choose a minimal representative geodesic loop $c\in C$. Let $\gamma \colon [0,\infty )\to M$ be a ray such that $\gamma (0)=c(0)$, and denote by $\theta$ the angle between $\dot{c}(0)$ and $\dot{\gamma}(0)$. Our first goal is to show that
\begin{equation}\label{ineq}
\theta \geq\theta _L:=
\arccos{(\cotanh{(L)}\tanh{(L/2)})}
\end{equation}
To see this, we work on the Riemannian universal cover $\tilde{M}$ of $M$ and lift $p$ to a point $\tilde{p}$. We also lift $c$ to a minimizing geodesic $\tilde{c}$ running from $\tilde{p}$ to $C\tilde{p}$ and $\gamma$ to a ray $\tilde{\gamma}$ starting at $\tilde{p}$. Note that the angle between $\dot{\tilde{c}}(0)$ and $\dot{\tilde{\gamma}}(0)$ is $\theta$. Moreover, we have
$$R=d_M(\gamma (R),p)\leq d_{\tilde{M}}(\tilde{\gamma}(R),C\tilde{p}).$$
Using this inequality and applying the Toponogov theorem to the hinge $(\tilde{p},\tilde{c}, \tilde{\gamma}\vert _{[0,R]})$, we get
$$\cosh{(R)}\leq \cosh ^2{(R)}-\sinh ^2{(R)}\cos{(\theta)}.$$
It follows that $\theta \geq \theta _R$. As the function $R\mapsto \theta _R$ is decreasing and $R\leq L$, we get (\ref{ineq}).

Now, we use an argument related to \cite[Lemma 2.7]{X2} and to the proof of \cite[Theorem 2]{dCX}. For each $u$ in the unit sphere $S^{n-1}\subset T_pM$, we denote by $\rho (u)$ the distance of $p$ to the cut point along the geodesic $t\mapsto \exp _p{(tu)}$. The function $\rho \colon S^{n-1}\to [0,\infty ]$ is continous. Furthermore, letting $du$ denote the induced measure on the unit sphere $S^{n-1}\subset T_pM$, we can write the volume form of the metric in geodesic polar coordinates around $p$ as $J(t,u)dt\,du$. 
Then the volume of a ball of radius $r$ around $p$ is
$$\vol {(B(p,r))}=\int _{S^{n-1}}\int _0^{\min{(r,\rho (u))}}J(t,u)\, dt\, du.$$
Using our assumption on the sectional curvature, we have $J(u,t)\leq \sinh ^{n-1}{(t)}$ by the Bishop-Gromov comparison theorem. Therefore we get
$$\vol {(B(p,r))}\leq \int _{S^{n-1}}\int _0^{\min{(r,\rho (u))}}\sinh ^{n-1}{(t)}\, dt\, du.$$
For $u\in S^{n-1}$, we denote by $\theta (u)$ the angle between $\dot{c}(0)$ and $u$. Let $\epsilon >0$ be a small fixed real number. Inequality (\ref{ineq}) tells us that if $\theta (u)\in[0,\theta _L-\epsilon]$, then $\rho (u)<\infty$. By continuity of $\rho$ and compacity of $\{ \theta (u)\leq \theta _L-\epsilon \}$, it follows that for some $\rho _0>0$, we have $\rho(u)\leq \rho _0$ for all $u$ such that $\theta (u)\leq \theta _L-\epsilon$. Hence, if $r$ is larger than $\rho _0$, we have
\begin{multline}
\vol {(B(p,r))}\leq \\
\int _{\{ \theta (u)\leq \theta _L-\epsilon\}}\int _0^{\rho _0}\sinh ^{n-1}{(t)}\, dt\, du + \int _{\{ \theta (u)\geq \theta _L-\epsilon\}}\int _0^{r}\sinh ^{n-1}{(t)}\, dt\, du.\nonumber
\end{multline}
Dividing by $V_{\H ^n}(r)=\omega _{n-1}\int _0^r \sinh ^{n-1}{(t)}\, dt$ and letting $r$ go to infinity, we get
$$\lim _{r\to \infty}{\frac{\vol{B(p,r)}}{V_{\H ^n}(r)}}\leq 1-\frac{V_{S^{n-1}}(\theta _L-\epsilon)}{\omega _{n-1}},$$
where $V_{S^{n-1}}(\theta _L-\epsilon)$ is the volume of a ball of radius $\theta _L-\epsilon$ in $S^{n-1}$. As $\epsilon$ is arbitrary, we can let it go to zero to achieve the proof.
\end{proof}

\Rk We have used the assumption $K\geq -1$ twice: the first time in order to apply the Toponogov comparison theorem on $\tilde{M}$ and the second time in order to get $J(u,t)\leq \sinh ^{n-1}(t)$ by using the Bishop-Gromov comparison theorem. Hence, if we happen to know a better bound on $J(u,t)$, we can improve a bit Theorem \ref{length}. Namely, let $f\colon (0,\infty)\to (0,\infty)$ be a measurable function and set $$F(r)=\int_0^r f(t)\, dt.$$ 
The proof of Theorem \ref{length} gives

\begin{prop}\label{pas}
Let $M^n$ be a complete $n-$dimensional non compact Riemannian manifold, and fix a point $p\in M$. Assume that the sectional curvature $K$ is bounded below $K\geq -1$, and that for all $t$ and $u$, we have $J(u,t)\leq f(t)$. For any real number $L>0$, there exists a constant $v=v(L,n)>0$ such that if
$$\limsup _{r\to \infty}{\frac{\vol{B(p,r)}}{F(r)}}\geq v$$
and $M$ is not simply connected, then the length of the shortest homotopically non trivial geodesic loop based at $p$ is bigger than $L$.
\end{prop}
For example, if the Ricci curvature is nonnegative, then the Bishop-Gromov theorem asserts that $J(u,t)\leq t^{n-1}$ and our proposition could be applied to manifolds of nonnegative Ricci curvature and sectional curvature bounded below. However, the constant $v$ of Proposition \ref{pas} can be easily seen to be bigger that $\omega_{n-1}/2$ and it follows from the work of Anderson  \cite{A} that a manifold of nonnegative Ricci curvature satisfying the volume growth of Proposition \ref{pas} is actually simply connected.  

\begin{proof}[Proof of Corollary \ref{Xia}]
We shall prove that $\pi _1(M,p)$ doesn't contain any  irreducible class and then apply Proposition \ref{generators}. If we choose the constant $v$ in the corollary to be the same as the one in Theorem \ref{length}, then we know that all irreducible elements must have length $\geq L$. If $C$ is such an element, with minimal representative geodesic loop $c$, then by Corollary \ref{critic}, the midpoint $c(L(C)/2)$ is a critical point of $d(p,.)$. However, it follows from the proof of \cite[Theorem 1]{X2} that if ii) and iii) are satisfied, then $d(p,.)$ doesn't have any critical point at distance $\geq L$.
\end{proof}


\section{Proofs of Theorems \ref{finite} and \ref{fg}}
\begin{proof}[Proof of Theorem \ref{finite}]
 We assume by contradiction that $M:=\tilde{M}/G$ is not compact and $G$ is infinite. Let $p\in M$ be the image of $\tilde{p}$ under the natural projection  $\tilde{M}\to M$. As $M$ is non compact, it is well known that there exists a ray $\gamma :[0,\infty )\to M$ emanating from $p$. We can lift $\gamma$ to a ray $\tilde{\gamma}$ in $\tilde{M}$ starting at $\tilde{p}$. Now if $\tilde{m}$ is a point in the orbit $G.\tilde{p}$ with $r:=d_{\tilde{M}}(\tilde{p},\tilde{m})$, we have
 $$d_{\tilde{M}}(\tilde{m},\tilde{\gamma}(r))\geq d_M(p,\gamma (r))=r.$$
But $G$ is assumed to be infinite, so that the orbit $G.\tilde{p}$ is not bounded. We can therefore find a sequence of points $\tilde{m_i}$ in $G.\tilde{p}$ such that $r_i=d_{\tilde{M}}(\tilde{p},\tilde{m_i})$ goes to infinity. This together with the inequality above contradicts the assumption of small diameter growth.
\end{proof}

\begin{proof}[Proof of Theorem \ref{fg}]
We assume on the contrary that $\pi _1(M,p)$ is not finitely generated. By Proposition \ref{generators}, we have an infinite sequence of irreducible elements $C_i$. Set $L_i=L(C_i)$; by Lemma \ref{fini} we know that $L_i$ goes to infinity as $i$ goes to infinity. For each $i$, let $c_i$ be a minimal representative geodesic loop in $C_i$. The vectors $\dot{c_i}(0)$ lie in the unit tangent sphere at $p$. As this sphere is compact, a subsequence, still to be denoted by $\dot{c_i}(0)$, converges to a unit tangent vector $u$ at $p$. By Proposition \ref{minimal} each $c_i$ is minimizing on $[0,L_i/2]$. Using this and the fact that $L_i$ goes to infinity, it is standard to conclude that $\gamma (t)= \exp _p{(tu)}$, $t\geq0$ is a ray in $M$ starting at $p$. We can lift $\gamma$ to a ray $\tilde{\gamma}$ in $\tilde{M}$ such that $\tilde{\gamma}(0)=\tilde{p}$. We can also lift each $c_i$ to a minimizing geodesic segment $\tilde{c_i}$ running from $\tilde{p}$ to $q_i:=C_i\tilde{p}$. Note that
$$d_{\tilde{M}}(q_i,\tilde{p})=L(C_i)=L_i.$$
Moreover, we have
$$d_{\tilde{M}}(q_i,\tilde{\gamma}(L_i))\geq d_M(p,\gamma (L_i))=L_i.$$
In other words, for all $i$ the hinge $(\tilde{p},\tilde{\gamma} |_{[0,L_i]},\tilde{c_i})$ is as described before Definition \ref{hinge}. However, its angle at $\tilde{p}$ is the same as the angle between $u$ and $\dot{c_i}(0)$, which by construction goes to zero if $i$ goes to infinity. Therefore $\tilde{M}$ does not have the thick hinges property at $\tilde{p}$, which is a contradiction.
\end{proof}



\begin{thebibliography}{99}

\bibitem[A]{A} M. Anderson, \emph{On the topology of complete manifolds of nonnegative Ricci curvature}, Topology 29 (1990),  no. 1, 41--55.

\bibitem[BW]{BW} I. Belegradek, G. Wei, \emph{Metrics of positive Ricci curvature on bundles}, Int. Math. Res. Not. (2004) no. 57, 3079--3096.

\bibitem[CW]{CW} J. Cao, J. Wang, \emph{private communication}.

\bibitem[Ch]{C} J. Cheeger, \emph{Critical points of distance functions and applications to geometry},  1--38, Lecture Notes in Math., 1504, Springer, Berlin, 1991. 

\bibitem[dCX]{dCX} M. do Carmo, C. Xia, \emph{Ricci curvature and the topology of open manifolds},  Math. Ann.  316  (2000),  no. 2, 391--400.

\bibitem[GP]{GP} R. Grimaldi, P. Pansu, \emph{Sur le degr\'e de diff\'erentiabilit\'e de la fonction croissance en dimension deux},  Boll. Un. Mat. Ital. B (7)  11  (1997),  no. 2, suppl., 25--38.

\bibitem[Gr1]{Gr} M. Gromov, \emph{Almost flat manifolds},  J. Differential Geom.  13  (1978), no. 2, 231--241.

\bibitem[Gr2]{Gr2} M. Gromov, \emph{Metric structures for Riemannian and non-Riemannian spaces}, Progress in Mathematics, 152, Birkh\"auser Boston, 1999.

\bibitem[Gro]{Gro} K. Grove, \emph{Critical point theory for distance functions},  Differential geometry: Riemannian geometry (Los Angeles, CA, 1990),  357--385, Proc. Sympos. Pure Math., 54, Part 3, Amer. Math. Soc., Providence, RI, 1993.

\bibitem[S]{Sh} Z. Shen, \emph{On complete manifolds of nonnegative $k$th-Ricci curvature}, Trans. Amer. Math. Soc.  338  (1993),  no. 1, 289--310.

\bibitem[So]{S} C. Sormani, \emph{Nonnegative Ricci curvature, small linear diameter growth and finite generation of fundamental groups},  J. Differential Geom.  54  (2000),  no. 3, 547--559. 


\bibitem[X]{X2} C. Xia, \emph{Complete manifolds with sectional curvature bounded below and large volume growth},  Bull. London Math. Soc.  34  (2002),  no. 2, 229--235.

\end{thebibliography}
\end{document}